\documentclass[11pt,reqno]{amsart}
\usepackage{amsmath,amssymb,amsxtra,latexsym,amsthm,amscd,amsfonts,mathrsfs,pb-diagram}
\usepackage[colorlinks=true,linkcolor=magenta,citecolor=blue]{hyperref}
\usepackage[margin=1in]{geometry}
\usepackage{enumitem} 
\usepackage{mathabx}
\usepackage{epsfig}

\newtheorem{dn}{Definition}[section]
\newtheorem{dl}{Theorem}[section]
\newtheorem{md}{Proposition}[section]
\newtheorem{bd}{Lemma}[section]
\newtheorem{hq}{Corollary}[section]
\newtheorem{nx}{Remark}[section]
\newtheorem{vd}{Example}[section]
\newcommand{\R}{\mathbb{R}}

\newcommand{\Z}{\mathbb{Z}}

\newcommand{\e}{\varepsilon}
\newcommand{\ity}{\infty}
\newcommand{\f}{\frac}

\newcommand{\bbd}{\begin{bd}}
\newcommand{\ebd}{\end{bd}}
\newcommand{\bdn}{\begin{dn}}
\newcommand{\edn}{\end{dn}}
\newcommand{\bhq}{\begin{hq}}
\newcommand{\ehq}{\end{hq}}
\newcommand{\bdl}{\begin{dl}}
\newcommand{\edl}{\end{dl}}
\newcommand{\bnx}{\begin{nx}}
\newcommand{\enx}{\end{nx}}
\newcommand{\bmd}{\begin{md}}
\newcommand{\emd}{\end{md}}
\newcommand{\bvd}{\begin{vd}}
\newcommand{\evd}{\end{vd}}

\title[Critical exponent for a $\sigma$-evolution system with frictional damping]{Critical exponent for a weakly coupled system of semi-linear $\sigma$-evolution equations with frictional damping}
\author{Tuan Anh Dao}
\address{Tuan Anh Dao \hfill\break
$\quad$ School of Applied Mathematics and Informatics, Hanoi University of Science and Technology, No.1 Dai Co Viet road, Hanoi, Vietnam \hfill\break
Faculty for Mathematics and Computer Science, TU Bergakademie Freiberg, Pr\"{u}ferstr. 9, 09596, Freiberg, Germany}
\email{anh.daotuan@hust.edu.vn}
\author{Trieu Duong Pham}
\address{Trieu Duong Pham \hfill\break
$\quad$ Department of Mathematics, Hanoi National University of Education, 136 Xuan Thuy, Hanoi, Vietnam}
\email{duongptmath@hnue.edu.vn}

\begin{document}
\subjclass{35B33, 35L56, 35S05}

\keywords{Frictional damping; $\sigma$-evolution equation; Weakly coupled system; Global existence; Critical exponent}
	
\begin{abstract}
We are interested in studying the Cauchy problem for a weakly coupled system of semi-linear $\sigma$-evolution equations with frictional damping. The main purpose of this paper is two-fold. We would like to not only prove the global (in time) existence of small data energy solutions but also indicate the blow-up result for Sobolev solutions when $\sigma$ is assumed to be any fractional number.
\end{abstract}

\maketitle

\section{Introduction} \label{Sec.main}
In this paper, let us consider the following Cauchy problem for weakly coupled system of semi-linear $\sigma$-evolution equations with frictional damping:
\begin{equation}
\begin{cases}
u_{tt}+ (-\Delta)^{\sigma_1} u+ u_t= |v|^p, &\quad x\in \R^n,\, t \ge 0, \\
v_{tt}+ (-\Delta)^{\sigma_2} v+ v_t= |u|^q, &\quad x\in \R^n,\, t \ge 0, \\
u(0,x)= u_0(x),\quad u_t(0,x)=u_1(x), &\quad x\in \R^n, \\
v(0,x)= v_0(x),\quad v_t(0,x)=v_1(x), &\quad x\in \R^n, \label{pt1.1}
\end{cases}
\end{equation}
for any $\sigma_1,\,\sigma_2\ge 1$ and for nonlinearities with powers $p,\, q >1$. We have in mind that the corresponding linear frictional damped $\sigma$-evolution equations of (\ref{pt1.1}) with vanishing right-hand side are in the following form:
\begin{equation}
w_{tt}+ (-\Delta)^{\sigma} w+ w_t=0,\quad w(0,x)= w_0(x),\quad w_t(0,x)= w_1(x), \label{pt1.2}
\end{equation}
with $\sigma= \sigma_1$ or $\sigma= \sigma_2$. One of the most typical important equations of (\ref{pt1.2}) with $\sigma=1$, well-known as the classical damped wave equation, has widely studied in numerous papers (see, for instance, \cite{Matsumura76,Matsumura77,Narazaki,TodorovaYordanov,Zhang} and the references therein). By the aid of the damping term $u_t$, the authors established decay estimates for solutions to the linear problem and applied these estimates to deal with the corresponding Cauchy problem for semi-linear equations with nonlinearity term $|u|^p$ or $|u|^{p-1}u$. Especially, related to nonlinearity term $|u|^p$ the global (in time) existence of energy solutions and a blow-up result were well-studied as well in the cited papers. After that, some previous results on the most well-known problem of a weakly coupled system (\ref{pt1.1}) with $\sigma_1=\sigma_2=1$, the so-called weakly coupled system of semi-linear classical damped wave equations, in the following form:
\begin{equation}
\begin{cases}
u_{tt}- \Delta u+ u_t= |v|^p,\quad  v_{tt}- \Delta v+ v_t= |u|^q, \\
u(0,x)= u_0(x),\quad u_t(0,x)=u_1(x),\quad v(0,x)= v_0(x),\quad v_t(0,x)=v_1(x).
\end{cases} \label{pt1.3}
\end{equation}
were explored in several papers (see \cite{Narazaki,NishiharaWakasugi,SunWang}). Concretely, Sun-Wang \cite{SunWang} and Narazaki \cite{Narazaki} proved the global (in time) existence of small data energy solutions to (\ref{pt1.3}) in low space dimensions $n=1,\,2,\,3$, provided that the following condition for the exponents $p,\,q$ is fulfilled:
$$ \frac{1+ \max\{p,\,q\}}{pq-1}< \frac{n}{2}. $$
Furthermore, if this condition is no longer true, then a result for nonexistence of global (in time) weak solutions to (\ref{pt1.3}) was indicated in \cite{Narazaki} and \cite{SunWang}. Afterwards, by using weighted energy estimates Nishihara-Wakasugi \cite{NishiharaWakasugi} extended these results for any space dimensions $n\ge 1$ with the same above condition.
\par As a natural extension of the classical damped wave equations, Takeda \cite{Takeda} considered (\ref{pt1.2}) when $\sigma$ is integer, the so-called polyharmonic damped wave equation. He has found a critical condition to ensure the global esxistence of small data solutions to a semi-linear Cauchy problem with the nonlinearity $|u|^p$ in the space dimension $n=1$. The limitation of the one-dimensional case is due to the technical difficulty. At present, the fact is that there seem not so many research papers concerning the study of (\ref{pt1.2}) for any fractional number $\sigma \ge 1$. Recently, Radu--Todorova-Yordanov \cite{RaduTodorovaYordanov} succeeded to derive some of sharp decay rates for solutions by proving a diffusion phenomenon in the following abstract setting:
\begin{equation}
w_{tt}+ B w+ w_t= 0, \label{pt1.4}
\end{equation}
where $B$ is a nonnegative self-adjoint operator. To establish these results, the authors have estimated, on the one hand, the decay rate of the difference between solutions to (\ref{pt1.4}) and those to the corresponding diffusion equation in terms of the flow $e^{-tB}$ of the operator $B$ by decomposing solutions to (\ref{pt1.4}) localized separately to low and high frequencies. On the other hand, they applied the Nash inequality and Markov property for the parabolic semigroup to yield explicit and sharp decay estimates. Quite recently, the application of these estimates to the operators $(-\Delta)^\sigma$ for any $\sigma \ge 1$ was investigated in the paper of Duong-Reissig \cite{DuongReissig}. The authors have discussed elementary results on the possible range of the admissible exponents in several of the corresponding semi-linear equations of (\ref{pt1.2}) with the nonlinearity term $\big||D|^a u\big|^p$.
\par To the best of author's knowledge, concerning the weakly coupled system of semi-linear frictional damped $\sigma$-evolution equations (\ref{pt1.1}) for any $\sigma\ge 1$, it seems that we still do not obtain any previous research manuscript. For this reason, the first main goal of this paper is to prove the global (in time) existence of small data energy solutions to (\ref{pt1.1}) for any $\sigma_1,\,\sigma_2 \ge 1$ by applying $(L^1 \cap L^2)- L^2$ estimates and $L^2- L^2$ estimates for solutions to (\ref{pt1.2}) from the recent paper of Duong-Reissig \cite{DuongReissig}. To do this, allowing loss of decay combined with using the fractional Gagliardo-Nirenberg inequality plays an important role in the treatment of the corresponding semi-linear equations. Moreover, in the present paper we want to explain that how the flexible choice of the parameters $\sigma_1,\,\sigma_2 \ge 1$ influences our global (in time) existence results and the range of admissible exponents $p,\,q$ as well. The fact is that taking into considerations the proof of blow-up results related to the fractional Laplacian operators $(-\Delta)^\sigma$, it seems difficult to directly apply the standard test function method (see, for example, \cite{Narazaki,NishiharaWakasugi,SunWang,TodorovaYordanov,Zhang}). Hence, our second main goal of this paper is to show a blow-up result to find the critical exponents by using a modified test function method (see more \cite{DaoReissig}) when $\sigma$ is assumed to be any fractional number.

\subsection{Notations}
Throughout this paper, we use the following notations.
\begin{itemize}[leftmargin=*]
\item We write $f\lesssim g$ when there exists a constant $C>0$ such that $f\le Cg$, and $f \approx g$ when $g\lesssim f\lesssim g$.
\item We denote $\widehat{f}(t,\xi):= \mathfrak{F}_{x\rightarrow \xi}\big(f(t,x)\big)$ as the Fourier transform with respect to the space variable of a function $f(t,x)$. As usual, the spaces $H^a$ and $\dot{H}^a$ with $a \ge 0$ stand for Bessel and Riesz potential spaces based on $L^2$ spaces. Here $\big<D\big>^{a}$ and $|D|^{a}$ denote the pseudo-differential operators with symbols $\big<\xi\big>^{a}$ and $|\xi|^{a}$, respectively.
\item For a given number $s \in \R$, we denote 
$$ [s]:= \max \big\{k \in \Z \,\, : \,\, k\le s \big\} \quad \text{ and }\quad [s]^+:= \max\{s,0\} $$
as its integer part and its positive part, respectively.
\item We put $\big< x\big>:= \sqrt{1+|x|^2}$, the so-called Japanese bracket of $x \in \R^n$.
\item Finally, we introduce the spaces
$\mathcal{A}^{\sigma}:= \big(L^1 \cap H^\sigma\big) \times \big(L^1 \cap L^2\big)$ with the norm
$$\|(u_0,u_1)\|_{\mathcal{A}^{\sigma}}:=\|u_0\|_{L^1}+ \|u_0\|_{H^\sigma}+ \|u_1\|_{L^1}+ \|u_1\|_{L^2}, \quad \text{ where }\sigma\ge 1. $$
\end{itemize}

\subsection{Main results}
At first, let us state the main results related to the global (in time) existence of small data energy solutions to (\ref{pt1.1}).
\bdl[\textbf{$\sigma_1 \ge \sigma_2$}] \label{dl1.1}
Let us assume $\sigma_1 \ge \sigma_2$. We assume that the conditions are satisfied
\begin{align}
&2 \le p,\, q < \ity & & & \text{ if }&\, n \le 2\sigma_2, \label{GN11A1} \\
&2 \le p \le \frac{n}{n- 2\sigma_2}, &\quad &2 \le q  < \ity & \text{ if }&\, 2\sigma_2 < n \le 2\sigma_1, \label{GN11A2} \\
&2 \le p \le \frac{n}{n- 2\sigma_2}, &\quad &2 \le q  \le \frac{n}{n- 2\sigma_1} & \text{ if }&\, 2\sigma_1 < n \le 4\sigma_2. \label{GN11A3}
\end{align}
Moreover, we suppose the following conditions:
\begin{equation} \label{exponent11A1}
\frac{1+ q}{(q-1)\big(\frac{\sigma_2}{\sigma_1}-1\big)+ pq-1}< \frac{n}{2\sigma_2},
\end{equation}
and
\begin{equation} \label{exponent11A2}
p \le 1+ \frac{2\sigma_2}{n} \le 1+ \frac{2\sigma_1}{n} < q.
\end{equation}
Then, there exists a constant $\e_0>0$ such that for any small data
$$ \big((u_0,u_1),\, (v_0,v_1) \big) \in \mathcal{A}^{\sigma_1} \times \mathcal{A}^{\sigma_2} \text{ satisfying the assumption } \|(u_0,u_1)\|_{\mathcal{A}^{\sigma_1}}+ \|(v_0,v_1)\|_{\mathcal{A}^{\sigma_2}} \le \e_0, $$
we have a uniquely determined global (in time) small data energy solution
$$ (u,v) \in \Big(C\big([0,\ity),H^{\sigma_1}\big)\cap C^1\big([0,\ity),L^2\big)\Big) \times \Big(C\big([0,\ity),H^{\sigma_2}\big)\cap C^1\big([0,\ity),L^2\big)\Big) $$
to (\ref{pt1.1}). The following estimates hold:
\begin{align}
\|u(t,\cdot)\|_{L^2} &\lesssim (1+t)^{-\frac{n}{4\sigma_1}+ \e(p,\sigma_2)} \big(\|(u_0,u_1)\|_{\mathcal{A}^{\sigma_1}}+ \|(v_0,v_1)\|_{\mathcal{A}^{\sigma_2}}\big), \label{decayrate11A1} \\
\big\||D|^{\sigma_1} u(t,\cdot)\big\|_{L^2} &\lesssim (1+t)^{-\frac{n}{4\sigma_1}-\frac{1}{2}+ \e(p,\sigma_2)} \big(\|(u_0,u_1)\|_{\mathcal{A}^{\sigma_1}}+ \|(v_0,v_1)\|_{\mathcal{A}^{\sigma_2}}\big), \label{decayrate11A2} \\
\|u_t(t,\cdot)\|_{L^2} &\lesssim (1+t)^{-\frac{n}{4\sigma_1}-1+ \e(p,\sigma_2)} \big(\|(u_0,u_1)\|_{\mathcal{A}^{\sigma_1}}+ \|(v_0,v_1)\|_{\mathcal{A}^{\sigma_2}}\big), \label{decayrate11A3} \\
\|v(t,\cdot)\|_{L^2} &\lesssim (1+t)^{-\frac{n}{4\sigma_2}} \big(\|(u_0,u_1)\|_{\mathcal{A}^{\sigma_1}}+ \|(v_0,v_1)\|_{\mathcal{A}^{\sigma_2}}\big), \label{decayrate11A4} \\
\big\||D|^{\sigma_2} v(t,\cdot)\big\|_{L^2} &\lesssim (1+t)^{-\frac{n}{4\sigma_2}-\frac{1}{2}} \big(\|(u_0,u_1)\|_{\mathcal{A}^{\sigma_1}}+ \|(v_0,v_1)\|_{\mathcal{A}^{\sigma_2}}\big), \label{decayrate11A5} \\
\|v_t(t,\cdot)\|_{L^2} &\lesssim (1+t)^{-\frac{n}{4\sigma_2}-1} \big(\|(u_0,u_1)\|_{\mathcal{A}^{\sigma_1}}+ \|(v_0,v_1)\|_{\mathcal{A}^{\sigma_2}}\big), \label{decayrate11A6}
\end{align}
where $\e(p,\sigma_2):= 1- \frac{n}{2\sigma_2}(p-1)+\e$ with a sufficiently small positive number $\e$.
\edl

\bdl[\textbf{$\sigma_2 \ge \sigma_1$}] \label{dl1.2}
Let us assume $\sigma_2 \ge \sigma_1$. We assume that the conditions are satisfied
\begin{align}
&2 \le p,\, q < \ity & & & \text{ if }&\, n \le 2\sigma_1, \label{GN12A1} \\
&2 \le p \le \frac{n}{n- 2\sigma_2}, &\quad &2 \le q  < \ity & \text{ if }&\, 2\sigma_1 < n \le 2\sigma_2, \label{GN12A2} \\
&2 \le p \le \frac{n}{n- 2\sigma_2}, &\quad &2 \le q  \le \frac{n}{n- 2\sigma_1} & \text{ if }&\, 2\sigma_2 < n \le 4\sigma_1. \label{GN12A3}
\end{align}
Moreover, we suppose the following conditions:
\begin{equation} \label{exponent12A1}
\frac{1+ p}{(p-1)\big(\frac{\sigma_1}{\sigma_2}-1\big)+ pq-1}< \frac{n}{2\sigma_1},
\end{equation}
and
\begin{equation} \label{exponent12A2}
q \le 1+ \frac{2\sigma_1}{n} \le 1+ \frac{2\sigma_2}{n} < p.
\end{equation}
Then, we have the same conclusions as in Theorem \ref{dl1.1}. But the estimates (\ref{decayrate11A1})-(\ref{decayrate11A6}) are modified in the following way:
\begin{align}
\|u(t,\cdot)\|_{L^2} &\lesssim (1+t)^{-\frac{n}{4\sigma_1}} \big(\|(u_0,u_1)\|_{\mathcal{A}^{\sigma_1}}+ \|(v_0,v_1)\|_{\mathcal{A}^{\sigma_2}}\big), \label{decayrate12A1} \\
\big\||D|^{\sigma_1} u(t,\cdot)\big\|_{L^2} &\lesssim (1+t)^{-\frac{n}{4\sigma_1}-\frac{1}{2}} \big(\|(u_0,u_1)\|_{\mathcal{A}^{\sigma_1}}+ \|(v_0,v_1)\|_{\mathcal{A}^{\sigma_2}}\big), \label{decayrate12A2} \\
\|u_t(t,\cdot)\|_{L^2} &\lesssim (1+t)^{-\frac{n}{4\sigma_1}-1} \big(\|(u_0,u_1)\|_{\mathcal{A}^{\sigma_1}}+ \|(v_0,v_1)\|_{\mathcal{A}^{\sigma_2}}\big), \label{decayrate12A3} \\
\|v(t,\cdot)\|_{L^2} &\lesssim (1+t)^{-\frac{n}{4\sigma_2}+ \e(q,\sigma_1)} \big(\|(u_0,u_1)\|_{\mathcal{A}^{\sigma_1}}+ \|(v_0,v_1)\|_{\mathcal{A}^{\sigma_2}}\big), \label{decayrate12A4} \\
\big\||D|^{\sigma_2} v(t,\cdot)\big\|_{L^2} &\lesssim (1+t)^{-\frac{n}{4\sigma_2}-\frac{1}{2}+ \e(q,\sigma_1)} \big(\|(u_0,u_1)\|_{\mathcal{A}^{\sigma_1}}+ \|(v_0,v_1)\|_{\mathcal{A}^{\sigma_2}}\big), \label{decayrate12A5} \\
\|v_t(t,\cdot)\|_{L^2} &\lesssim (1+t)^{-\frac{n}{4\sigma_2}-1+ \e(q,\sigma_1)} \big(\|(u_0,u_1)\|_{\mathcal{A}^{\sigma_1}}+ \|(v_0,v_1)\|_{\mathcal{A}^{\sigma_2}}\big), \label{decayrate12A6}
\end{align}
where $\e(q,\sigma_1):= 1- \frac{n}{2\sigma_1}(q-1)+\e$ with a sufficiently small positive number $\e$.
\edl

\begin{nx}
\fontshape{n}
\selectfont
It is clear that because of the conditions (\ref{exponent11A1}) and (\ref{exponent12A1}), both $\e(p,\sigma_2)$ and $\e(q,\sigma_1)$ appearing in Theorems \ref{dl1.1} and \ref{dl1.2} are non-negative.
\end{nx}

\begin{nx}
\fontshape{n}
\selectfont
Let us explain the interaction between the flexible choice of the parameters $\sigma_1,\,\sigma_2$ and the admissible exponents $p,\,q$ in Theorems \ref{dl1.1} and \ref{dl1.2}. It is obvious that Theorem \ref{dl1.1} concerns the case $p< q$ corresponding to the case $\sigma_1\ge \sigma_2$, whereas the inverse case $q< p$ is of interest in Theorem \ref{dl1.2} corresponding to the case $\sigma_2\ge \sigma_1$. Here we can see that the different choice of the parameters $\sigma_1,\,\sigma_2$ affects our admissible exponents $p,\,q$ remarkably. Besides, taking account of the special case of parameters $\sigma_1= \sigma_2=: \sigma$ we may observe that the conditions from (\ref{GN11A1}) to (\ref{GN11A3}) are the same as those from (\ref{GN12A1}) to (\ref{GN12A3}). Hence, it is reasonable to re-write the assumptions (\ref{exponent11A1}) and (\ref{exponent11A2}) in Theorem \ref{dl1.1}, (\ref{exponent12A1}) and (\ref{exponent12A2}) in Theorem \ref{dl1.2} in the following common form:
$$ \frac{1+ \max\{p,\,q\}}{pq-1}< \frac{n}{2\sigma}, $$
and
$$ \min\{p,\,q\} \le 1+ \frac{2\sigma}{n} < \max\{p,\,q\}. $$
\end{nx}

\begin{nx}
\fontshape{n}
\selectfont
Here we want to underline that we have derived the global (in time) existence of small data energy solutions to (\ref{pt1.1}) in this paper. For the purpose of further considerations, we may expect to obtain several results concerning Sobolev solutions, energy solutions with a higher regularity or large regular solutions to (\ref{pt1.1}) as we have proved in the paper \cite{Dao} by applying some new tools from Harmonic Analysis such as the fractional Leibniz rule, the fractional chain rule, the fractional powers rule and the fractional Sobolev embedding.
\end{nx}

Finally, the third result is concerned with the following blow-up result to indicate the sharpness of our exponents to (\ref{pt1.1}).

\bdl[\textbf{Blow-up}] \label{dl1.3}
Let $\sigma_1= \sigma_2=: \sigma \ge 1$ be a fractional number. We assume that we choose the initial data $u_0=v_0=0$ and $u_1,\,v_1 \in L^1$ satisfying the following relations:
\begin{equation} \label{optimal13.1}
\int_{\R^n} u_1(x)dx > \epsilon_1 \quad \text{ and }\quad \int_{\R^n} v_1(x)dx > \epsilon_2,
\end{equation}
where $\epsilon_1$ and $\epsilon_2$ are suitable nonnegative constants. Moreover, we suppose the following condition:
\begin{equation}
\frac{n}{2\sigma} \le \frac{1+ \max\{p,\,q\}}{pq-1} \label{optimal13.2}.
\end{equation}
Then, there is no global (in time) Sobolev solution $(u,v) \in C\big([0,\infty),L^2\big) \times C\big([0,\infty),L^2\big)$ to (\ref{pt1.1}).
\edl

\begin{nx}
\fontshape{n}
\selectfont
We can see that if we choose $\sigma_1= \sigma_2= \sigma$ in Theorems \ref{dl1.1} and \ref{dl1.2}, then from Theorem \ref{dl1.3} it follows that the exponents $p,\,q$ given by
$$ \frac{1+ \max\{p,\,q\}}{pq-1}= \frac{n}{2\sigma} $$
are really critical. Especially, by choosing $\sigma_1= \sigma_2=1$ in Theorems \ref{dl1.1}, \ref{dl1.2} and \ref{dl1.3} we recognize that our obtained results are consistent with those in \cite{Narazaki,NishiharaWakasugi,SunWang}.
\end{nx}

\textbf{The structure of this article is organized as follows}: In Section \ref{Sec.Pre}, we collect $(L^1 \cap L^2)- L^2$ estimates and $L^2- L^2$ estimates for solutions to (\ref{pt1.2}) from the recent paper of Duong-Reissig \cite{DuongReissig}, and give some necessary properties of the modified test function method as well. We present the proofs of our global (in time) existence results to (\ref{pt1.1}) in Section \ref{Sec.Global existence}. Finally, Section \ref{Sec.Blow-up} is devoted to prove the blow-up result.

\section{Preliminaries} \label{Sec.Pre}

\subsection{Linear estimates}
Main goal of this section is to collect $(L^1 \cap L^2)- L^2$ and $L^2- L^2$ estimates for solutions and some of their derivatives to (\ref{pt1.2}) from the recent paper of Duong-Reissig \cite{DuongReissig}. At first, applying partial Fourier transformation to (\ref{pt1.2}) we derive the following Cauchy problem:
\begin{equation}
\widehat{w}_{tt}+ \widehat{w}_t+ |\xi|^{2\sigma} \widehat{w}=0,\quad \widehat{w}(0,\xi)= \widehat{w}_0(\xi),\quad \widehat{w}_t(0,\xi)= \widehat{w}_1(\xi). \label{pt2.1}
\end{equation}
The characteristic roots are
$$ \lambda_{1,2}=\lambda_{1,2}(\xi)= \f{1}{2}\Big(-1 \pm \sqrt{1-4|\xi|^{2\sigma}}\Big). $$
The solutions to (\ref{pt2.1}) are presented by the following formula (here we assume $\lambda_{1}\neq \lambda_{2}$):
\begin{align*}
\widehat{w}(t,\xi)&= \frac{\lambda_1 e^{\lambda_2 t}-\lambda_2 e^{\lambda_1 t}}{\lambda_1- \lambda_2}\widehat{w}_0(\xi)+ \frac{e^{\lambda_1 t}-e^{\lambda_2 t}}{\lambda_1- \lambda_2}\widehat{w}_1(\xi) \\
&=: \widehat{K}_{0,\sigma}(t,\xi)\widehat{w}_0(\xi)+ \widehat{K}_{1,\sigma}(t,\xi)\widehat{w}_1(\xi),
\end{align*}
that is, we may write the solutions to (\ref{pt1.2}) in the form
$$ w(t,x)=K_{0,\sigma}(t,x) \ast_{x} w_0(x)+ K_{1,\sigma}(t,x) \ast_{x} w_1(x), $$
where
$$K_{0,\sigma}(t,x)= \mathfrak{F}^{-1}_{\xi \to x}\big(\widehat{K}_{0,\sigma}(t,\xi)\big) \quad \text{ and }\quad K_{1,\sigma}(t,x)= \mathfrak{F}^{-1}_{\xi \to x}\big(\widehat{K}_{1,\sigma}(t,\xi)\big). $$
From the cited paper we may conclude the following statements.
\bmd[Proposition 2.1 in \cite{DuongReissig}]  \label{md2.1}
Let $\sigma= \sigma_k \ge 1$ with $k=1,\,2$ in (\ref{pt1.2}). The solutions to (\ref{pt1.2}) satisfy the $(L^1 \cap L^2)-L^2$ estimates
\begin{align*}
\|w(t,\cdot)\|_{L^2} &\lesssim (1+t)^{-\frac{n}{4\sigma_k}}\|w_0\|_{L^1 \cap L^2}+ (1+t)^{-\frac{n}{4\sigma_k}}\|w_1\|_{L^1 \cap H^{-{\sigma_k}}}, \\
\big\||D|^{\sigma_k} w(t,\cdot)\big\|_{L^2} &\lesssim (1+t)^{-\frac{n}{4\sigma_k}- \frac{1}{2}}\|w_0\|_{L^1 \cap H^{\sigma_k}}+ (1+t)^{-\frac{n}{4\sigma_k}- \frac{1}{2}}\|w_1\|_{L^1 \cap L^2}, \\
\|w(t,\cdot)\|_{L^2} &\lesssim (1+t)^{-\frac{n}{4\sigma_k}-1}\|w_0\|_{L^1 \cap H^{\sigma_k}}+ (1+t)^{-\frac{n}{4\sigma_k}-1}\|w_1\|_{L^1 \cap L^2},
\end{align*}
and the $L^2-L^2$ estimates
\begin{align*}
\|w(t,\cdot)\|_{L^2} &\lesssim \|w_0\|_{L^2}+ (1+t)\|w_1\|_{L^2}, \\
\big\||D|^{\sigma_k} w(t,\cdot)\big\|_{L^2} &\lesssim (1+t)^{-\frac{1}{2}}\|w_0\|_{H^{\sigma_k}}+ (1+t)^{-\frac{1}{2}}\|w_1\|_{L^2}, \\
\|w(t,\cdot)\|_{L^2} &\lesssim (1+t)^{-1}\|w_0\|_{H^{\sigma_k}}+  \|w_1\|_{L^2},
\end{align*}
for all space dimensions $n \ge 1$.
\emd

\begin{nx}
\fontshape{n}
\selectfont
The statements in Propositions \ref{md2.1} are key tools to prove the global (in time) existence results for (\ref{pt1.1}) in Section \ref{Sec.Global existence}.
\end{nx}

\subsection{A modified test function}
In this section, we would like to present some auxiliary properties of the modified test function $\psi= \psi(x):= \big< x\big>^{-r}$ for some $r>0$ from the recent paper of Dao-Reissig \cite{DaoReissig} which plays an essential role to prove our blow-up result in Section \ref{Sec.Blow-up}.

\bbd[Lemma 2.3 in \cite{DaoReissig}] \label{lemma2.1}
Let $\gamma \ge 1$ be a fractional number and $s:= \gamma- [\gamma]$. Let $r>0$. Then, the following estimates hold for all $x \in \R^n$:
$$ \big|(-\Delta)^\gamma \big< x\big>^{-r}\big| \lesssim
\begin{cases}
\big< x\big>^{-r-2\gamma} &\quad \text{ if }\quad 0< r+2[\gamma]< n, \\
\big< x\big>^{-n-2s}\log(e+|x|) &\quad \text{ if }\quad r+2[\gamma]= n, \\
\big< x\big>^{-n-2s} &\quad \text{ if }\quad r+2[\gamma]> n.
\end{cases} $$
\ebd

\bbd \label{lemma2.2}
Let $\gamma \ge 1$ be a fractional number. Let $\psi:= \psi(x)= \big< x\big>^{-r}$ for some $r>\frac{n}{2}$. For any $R>0$, let $\psi_R$ be a function defined by
$$ \psi_R(x):= \psi\big(R^{-1} x\big)\quad \text{ for all }x \in \R^n. $$
Then, $(-\Delta)^\gamma (\psi_R)$ satisfies the following scaling properties for all $x \in \R^n$:
\begin{equation*}
(-\Delta)^\gamma (\psi_R)(x)= R^{-2\gamma} \big((-\Delta)^\gamma \psi \big)\big(R^{-1} x\big).
\end{equation*}
\ebd
\begin{proof}
At first, a direct computation leads to the following formula:
\begin{equation}
-\Delta \big< x\big>^{-\ell}= \ell\Big((n-\ell-2)\big< x\big>^{-\ell-2} + (\ell+2)\big< x\big>^{-\ell-4}\Big) \quad \text{ for any } \ell>0. \label{RepresentationFormula^1}
\end{equation}
After performing $[\gamma]$ steps of (\ref{RepresentationFormula^1}) and using induction argument, we arrive at the following formula:
\begin{align}
(-\Delta)^{[\gamma]} \big< x\big>^{-r}&= (-1)^{[\gamma]} \prod_{j=0}^{[\gamma]-1}(r+2j)\Big(\prod_{j=1}^{[\gamma]}(-n+r+2j)\big< x\big>^{-r-2[\gamma]} \nonumber \\
&\hspace{4cm}- C^1_{[\gamma]} \prod_{j=2}^{[\gamma]}(-n+r+2j)(r+2[\gamma])\big< x\big>^{-r-2[\gamma]-2} \nonumber \\
&\hspace{4cm}+ C^2_{[\gamma]} \prod_{j=3}^{[\gamma]}(-n+r+2j)(r+2m)(r+2[\gamma]+2)\big< x\big>^{-r-2[\gamma]-4} \nonumber \\
&\hspace{4cm}+\cdots+ (-1)^{[\gamma]} \prod_{j=0}^{[\gamma]-1}(r+2[\gamma]+2j)\big< x\big>^{-r-4[\gamma]}\Big). \label{RepresentationFormula^m}
\end{align}
Carrying out the change of variables $\tilde{x}:= R^{-1}x$ we obtain
$$ (-\Delta)^{[\gamma]}\psi_R(x)= R^{-2[\gamma]}(-\Delta)^{[\gamma]}(\psi)(\tilde{x}) $$
because $[\gamma]$ is an integer number. Employing the formula (\ref{RepresentationFormula^m}) we may re-write as follows:
\begin{align*}
(-\Delta)^{[\gamma]}\psi_R(x)&= (-1)^{[\gamma]} R^{-2[\gamma]} \prod_{j=0}^{[\gamma]-1}(r+2j)\Big(\prod_{j=1}^{[\gamma]}(-n+r+2j)\big< \tilde{x}\big>^{-r-2[\gamma]} \\
&\hspace{3cm}- C^1_{[\gamma]} \prod_{j=2}^{[\gamma]}(-n+r+2j)(r+2[\gamma])\big< \tilde{x}\big>^{-r-2[\gamma]-2} \\
&\hspace{3cm}+ C^2_{[\gamma]} \prod_{j=3}^{[\gamma]}(-n+r+2j)(r+2[\gamma])(r+2[\gamma]+2)\big< \tilde{x}\big>^{-r-2[\gamma]-4} \\
&\hspace{3cm}+\cdots+ (-1)^{[\gamma]} \prod_{j=0}^{[\gamma]-1}(r+2[\gamma]+2j)\big< \tilde{x}\big>^{-r-4[\gamma]}\Big).
\end{align*}
Let us introduce the following auxiliary functions:
$$ \psi_k(x):= \big< x\big>^{-r-2[\gamma]-2k}\quad \text{ and }\quad \psi_{k,R}(x):= \varphi_k(R^{-1}x)=\big< \tilde{x}\big>^{-r-2[\gamma]-2k} $$
with $k=0,\cdots,[\gamma]$. Now we write $\gamma= [\gamma]+ s$, where $s:= \gamma- [\gamma] \in (0,1)$ since $\gamma$ is fractional. Using the relation
$$ (-\Delta)^{\gamma}\psi_R(x)= (-\Delta)^s \big((-\Delta)^{[\gamma]}\psi_R(x)\big) $$
we get
\begin{align*}
(-\Delta)^{\gamma}\psi_R(x)&= (-1)^{[\gamma]} R^{-2[\gamma]} \prod_{j=0}^{[\gamma]-1}(r+2j)\Big(\prod_{j=1}^{[\gamma]}(-n+r+2j)\, (-\Delta)^s (\psi_{0,R})(x) \\
&\hspace{3cm}- C^1_{[\gamma]} \prod_{j=2}^{[\gamma]}(-n+r+2j)(r+2[\gamma])\, (-\Delta)^s (\psi_{1,R})(x) \\
&\hspace{3cm}+ C^2_{[\gamma]} \prod_{j=3}^{[\gamma]}(-n+r+2j)(r+2[\gamma])(r+2[\gamma]+2)\, (-\Delta)^s (\psi_{2,R})(x) \\
&\hspace{3cm}+\cdots+ (-1)^{[\gamma]} \prod_{j=0}^{[\gamma]-1}(r+2[\gamma]+2j)\, (-\Delta)^s (\psi_{[\gamma],R})(x)\Big).
\end{align*}
Consequently, applying Lemma \ref{UsefulLemma} we may conclude
\begin{align*}
(-\Delta)^{\gamma}\psi_R(x)&= (-1)^{[\gamma]} R^{-2[\gamma]-2s} \prod_{j=0}^{[\gamma]-1}(r+2j)\Big(\prod_{j=1}^{[\gamma]}(-n+r+2j)\, (-\Delta)^s (\psi_0)(\tilde{x}) \\
&\hspace{3cm}- C^1_{[\gamma]} \prod_{j=2}^{[\gamma]}(-n+r+2j)(r+2[\gamma])\, (-\Delta)^s (\psi_1)(\tilde{x}) \\
&\hspace{3cm}+ C^2_{[\gamma]} \prod_{j=3}^{[\gamma]}(-n+r+2j)(q+2[\gamma])(r+2[\gamma]+2)\, (-\Delta)^s (\psi_2)(\tilde{x}) \\
&\hspace{3cm}+\cdots+ (-1)^{[\gamma]} \prod_{j=0}^{[\gamma]-1}(r+2[\gamma]+2j)\, (-\Delta)^s (\psi_{[\gamma]})(\tilde{x})\Big) \\
&= R^{-2\gamma} (-\Delta)^{\gamma}(\psi)(\tilde{x}).
\end{align*}
Therefore, this completes our proof.
\end{proof}

\bbd[Lemma 2.7 in \cite{DaoReissig}] \label{lemma2.3}
Let $s \in \R$. Let $\psi_1= \psi_1(x) \in H^s$ and $\psi_2= \psi_2(x) \in H^{-s}$. Then, the following relation holds:
$$ \int_{\R^n}\psi_1(x)\,\psi_2(x)dx= \int_{\R^n}\widehat{\psi}_1(\xi)\,\widehat{\psi}_2(\xi)d\xi. $$
\ebd

\section{Proof of global existence results} \label{Sec.Global existence}

\subsection{Proof of Theorem \ref{dl1.1}}
First, we choose the data spaces $(u_0,u_1) \in \mathcal{A}^{\sigma_1}$ and $(v_0,v_1) \in \mathcal{A}^{\sigma_2}$. We introduce the family $\{X(t)\}_{t>0}$ of the solution spaces
$$ X(t):= \Big(C\big([0,t],H^{\sigma_1}\big)\cap C^1\big([0,t],L^2\big)\Big) \times \Big(C\big([0,t],H^{\sigma_2}\big)\cap C^1\big([0,t],L^2\big)\Big) $$
with the norm
\begin{align*}
\|(u,v)\|_{X(t)}:= \sup_{0\le \tau \le t} \Big( &f_1(\tau)^{-1}\|u(\tau,\cdot)\|_{L^2} + f_2(\tau)^{-1}\big\||D|^{\sigma_1} u(\tau,\cdot)\big\|_{L^2}+ f_3(\tau)^{-1}\|u_t(\tau,\cdot)\|_{L^2} \\
&+g_1(\tau)^{-1}\|v(\tau,\cdot)\|_{L^2} + g_2(\tau)^{-1}\big\||D|^{\sigma_2} v(\tau,\cdot)\big\|_{L^2}+ g_3(\tau)^{-1}\|v_t(\tau,\cdot)\|_{L^2} \Big),
\end{align*}
where
\begin{align}
&f_1(\tau)= (1+\tau)^{-\frac{n}{4\sigma_1}+ \e(p,\sigma_2)},\,\,\, f_2(\tau)= (1+\tau)^{-\frac{n}{4\sigma_1}-\frac{1}{2}+ \e(p,\sigma_2)}, \,\,\, f_3(\tau)=(1+\tau)^{-\frac{n}{4\sigma_1}-1+ \e(p,\sigma_2)}, \label{pt4.1} \\
&g_1(\tau)= (1+\tau)^{-\frac{n}{4\sigma_2}},\quad g_2(\tau)= (1+\tau)^{-\frac{n}{4\sigma_2}-\frac{1}{2}}, \quad g_3(\tau)=(1+\tau)^{-\frac{n}{4\sigma_2}-1}. \label{pt4.2}
\end{align}
By recalling the fundamental solutions $K_{0,\sigma}(t,x)$ and $K_{1,\sigma}(t,x)$ as defined in Section \ref{Sec.Pre}, we write the solutions of the corresponding linear Cauchy problems with vanishing right-hand sides to (\ref{pt1.1}) in the following form:
$$\begin{cases}
u^{ln}(t,x)=K_{0,\sigma_1}(t,x) \ast_{x} u_0(x)+ K_{1,\sigma_1}(t,x) \ast_{x} u_1(x), \\ 
v^{ln}(t,x)=K_{0,\sigma_2}(t,x) \ast_{x} v_0(x)+ K_{1,\sigma_2}(t,x) \ast_{x} v_1(x).
\end{cases}$$
Since we are interested in dealing with the semi-linear models with constant coefficients in the linear part, we apply Duhamel's principle to obtain the following formal implicit representation of solutions to (\ref{pt1.1}):
$$\begin{cases}
u(t,x)= u^{ln}(t,x) + \displaystyle\int_0^t K_{1,\sigma_1}(t-\tau,x) \ast_x |v(\tau,x)|^p d\tau=: u^{ln}(t,x)+ u^{nl}(t,x), \\ 
v(t,x)= v^{ln}(t,x) + \displaystyle\int_0^t K_{1,\sigma_2}(t-\tau,x) \ast_x |u(\tau,x)|^q d\tau=: v^{ln}(t,x)+ v^{nl}(t,x).
\end{cases}$$
We define the following operator for all $t>0$:
\begin{align*}
N &: \quad X(t) \longrightarrow X(t) \\ 
N(u,v)(t,x) &= \big(u^{ln}(t,x)+ u^{nl}(t,x), v^{ln}(t,x)+ v^{nl}(t,x)\big).
\end{align*}
We shall indicate that the operator $N$ fulfills the following two inequalities:
\begin{align}
\|N(u,v)\|_{X(t)} &\lesssim \|(u_0,u_1)\|_{\mathcal{A}^{\sigma_1}}+ \|(v_0,v_1)\|_{\mathcal{A}^{\sigma_2}}+ \|(u,v)\|^p_{X(t)}+ \|(u,v)\|^q_{X(t)}, \label{pt4.3} \\
\|N(u,v)-N(\bar{u},\bar{v})\|_{X(t)} &\lesssim \|(u,v)-(\bar{u},\bar{v})\|_{X(t)} \Big(\|(u,v)\|^{p-1}_{X(t)}+ \|(\bar{u},\bar{v})\|^{p-1}_{X(t)} \nonumber \\
&\hspace{6cm}+ \|(u,v)\|^{q-1}_{X(t)}+ \|(\bar{u},\bar{v})\|^{q-1}_{X(t)}\Big). \label{pt4.4}
\end{align}
Then, we may conclude global (in time) existence results of small data solutions by applying Banach's fixed point theorem.
\par In the first step, from the statements in Proposition \ref{md2.1} and the definition of the norm in $X(t)$  we get
\begin{equation*}
\big\|(u^{ln}, v^{ln})\big\|_{X(t)} \lesssim \|(u_0,u_1)\|_{\mathcal{A}^{\sigma_1}}+ \|(v_0,v_1)\|_{\mathcal{A}^{\sigma_2}}.
\end{equation*}
For this reason, to prove (\ref{pt4.3}) it is sufficient to show the following inequality:
\begin{equation} \label{pt4.5}
\big\|(u^{nl}, v^{nl})\big\|_{X(t)} \lesssim \|(u,v)\|^p_{X(t)}+ \|(u,v)\|^q_{X(t)}.
\end{equation}
\par Let us prove the inequality (\ref{pt4.5}). In order to control $u^{nl}$, we use the $(L^1 \cap L^2)- L^2$ estimates from Proposition \ref{md2.1} to derive the following estimate:
$$ \big\|u^{nl}(t,\cdot)\big\|_{L^2}\lesssim \int_0^t (1+t-\tau)^{-\frac{n}{4\sigma_1}}\big\||v(\tau,\cdot)|^p\big\|_{L^1 \cap L^2}d\tau. $$
Therefore, it is necessary to require the estimates for $|v(\tau,x)|^p$ in $L^1$ and $L^2$ as follows:
$$ \big\||v(\tau,\cdot)|^p\big\|_{L^1}= \|v(\tau,\cdot)\|^p_{L^p} \quad \text{ and }\quad \big\||v(\tau,\cdot)|^p\big\|_{L^2}= \|v(\tau,\cdot)\|^p_{L^{2p}}. $$
Applying the fractional Gagliardo-Nirenberg inequality from Proposition \ref{fractionalGagliardoNirenberg} gives
\begin{align}
\big\||v(\tau,\cdot)|^p\big\|_{L^1} &\lesssim (1+\tau)^{-\frac{n}{2\sigma_2}(p-1)}\|(u,v)\|^p_{X(\tau)}, \label{t11.1} \\
\big\||v(\tau,\cdot)|^p\big\|_{L^2} &\lesssim (1+\tau)^{-\frac{n}{2\sigma_2}(p-\frac{1}2)}\|(u,v)\|^p_{X(\tau)}, \label{t11.2}
\end{align}
provided that the conditions from (\ref{GN11A1}) to (\ref{GN11A3}) are satisfied for $p$. From the both above estimates, we have
\begin{align*}
\big\|u^{nl}(t,\cdot)\big\|_{L^2} &\lesssim (1+t)^{-\frac{n}{4\sigma_1}}\|(u,v)\|^p_{X(t)} \int_0^{t/2}(1+\tau)^{-\frac{n}{2\sigma_2}(p-1)} d\tau \\
&\qquad + (1+t)^{-\frac{n}{2\sigma_2}(p-1)}\|(u,v)\|^p_{X(t)} \int_{t/2}^t (1+t-\tau)^{-\frac{n}{4\sigma_1}}d\tau.
\end{align*}
Here we used the relation
\begin{equation}
(1+t-\tau) \approx (1+t) \text{ if }\tau \in [0,t/2], \text{ and } (1+\tau) \approx (1+t) \text{ if }\tau \in [t/2,t]. \label{t11.3}
\end{equation}
Because of the condition $p \le 1+ \frac{2\sigma_2}{n}$ in (\ref{exponent11A2}), it is clear that the term $(1+\tau)^{-\frac{n}{2\sigma_2}(p-1)}$ is not integrable. Consequently, we obtain
\begin{align*}
(1+t)^{-\frac{n}{4\sigma_1}} \int_0^{t/2}(1+\tau)^{-\frac{n}{2\sigma_2}(p-1)} d\tau &\lesssim
\begin{cases}
(1+t)^{-\frac{n}{4\sigma_1}+1-\frac{n}{2\sigma_2}(p-1)} &\text{ if }p< 1+ \frac{2\sigma_2}{n} \\
(1+t)^{-\frac{n}{4\sigma_1}+\e} &\text{ if }p= 1+ \frac{2\sigma_2}{n}
\end{cases} \\ 
&\lesssim (1+t)^{-\frac{n}{4\sigma_1}+ \e(p,\sigma_2)},
\end{align*}
where $\e$ is a sufficiently small positive number. Furthermore, we can see that $\frac{n}{4\sigma_1}\le 1$ due to the condition $n\le 4\sigma_2$ in (\ref{GN11A3}) and the assumption $\sigma_1 \ge \sigma_2$. Thus, it follows that
\begin{align*}
(1+t)^{-\frac{n}{2\sigma_2}(p-1)} \int_{t/2}^t (1+t-\tau)^{-\frac{n}{4\sigma_1}}d\tau &\lesssim
\begin{cases}
(1+t)^{-\frac{n}{2\sigma_2}(p-1)+1-\frac{n}{4\sigma_1}} &\text{ if } n< 4\sigma_1 \\
(1+t)^{-\frac{n}{2\sigma_2}(p-1)+\e} &\text{ if } n= 4\sigma_1
\end{cases} \\
&\lesssim (1+t)^{-\frac{n}{4\sigma_1}+ \e(p,\sigma_2)}
\end{align*}
with a sufficiently small positive number $\e$. From the both above estimates, we may arrive at the following estimate:
$$\big\|u^{nl}(t,\cdot)\big\|_{L^2} \lesssim (1+t)^{-\frac{n}{4\sigma_1}+ \e(p,\sigma_2)} \|(u,v)\|^p_{X(t)}. $$
To deal with $|D|^{\sigma_1}u^{nl}$, we use the $(L^1 \cap L^2)- L^2$ estimates if $\tau \in [0,t/2]$ and the $L^2-L^2$ estimates if $\tau \in [t/2,t]$ from Proposition \ref{md2.1} to get the following estimate:
\begin{align*}
\big\||D|^{\sigma_1}u^{nl}(t,\cdot)\big\|_{L^2}&\lesssim \int_0^{t/2} (1+t-\tau)^{-\frac{n}{4\sigma_1}-\frac{1}{2}}\big\||v(\tau,\cdot)|^p\big\|_{L^1 \cap L^2}d\tau+\int_{t/2}^t (1+t-\tau)^{-\frac{1}{2}}\big\||v(\tau,\cdot)|^p\big\|_{L^2}d\tau \\
&\lesssim (1+t)^{-\frac{n}{4\sigma_1}-\frac{1}{2}}\|(u,v)\|^p_{X(t)} \int_0^{t/2}(1+\tau)^{-\frac{n}{2\sigma_2}(p-1)}d\tau \\
&\qquad +(1+t)^{-\frac{n}{2\sigma_2}(p-\frac{1}2)}\|(u,v)\|^p_{X(t)} \int_{t/2}^t (1+t-\tau)^{-\frac{1}{2}}d\tau.
\end{align*}
Here we used again the estimates (\ref{t11.1}) and (\ref{t11.2}) combined with the relation (\ref{t11.3}). The first integral will be handled as we did to estimate $u^{nl}$. As a result, we derive
$$ (1+t)^{-\frac{n}{4\sigma_1}-\frac{1}{2}} \int_0^{t/2}(1+\tau)^{-\frac{n}{2\sigma_2}(p-1)} d\tau \lesssim (1+t)^{-\frac{n}{4\sigma_1}-\frac{1}{2}+\e(p_1,\sigma_2)}. $$
For the second integral, we can proceed as follows:
$$ (1+t)^{-\frac{n}{2\sigma_2}(p-\frac{1}2)} \int_{t/2}^t (1+t-\tau)^{-\frac{1}{2}}d\tau \lesssim (1+t)^{-\frac{n}{2\sigma_2}(p-\frac{1}2)+\frac{1}{2}}\lesssim (1+t)^{-\frac{n}{4\sigma_1}-\frac{1}{2}+\e(p,\sigma_2)} $$
since $\sigma_1 \ge \sigma_2$. Therefore, we may conclude the following estimate:
$$\big\||D|^{\sigma_1}u^{nl}(t,\cdot)\big\|_{L^2} \lesssim (1+t)^{-\frac{n}{4\sigma_1}-\frac{1}{2}+\e(p,\sigma_2)} \|(u,v)\|^p_{X(t)}. $$
Similar to the treatment of $|D|^{\sigma_1}u^{nl}$, we also obtain
$$\big\|u_t^{nl}(t,\cdot)\big\|_{L^2} \lesssim (1+t)^{-\frac{n}{4\sigma_1}-1+\e(p,\sigma_2)}\|(u,v)\|^p_{X(t)}. $$
Analogously, we may arrive at the following estimates for $j,k=0,1$ with $(j,k)\ne (1,1)$:
$$\big\|\partial^j_t |D|^{k\sigma_2}v^{nl}(t,\cdot)\big\|_{L^2} \lesssim (1+t)^{-\frac{n}{4\sigma_2}- \frac{k}{2}-j}\|(u,v)\|^q_{X(t)}, $$
where the conditions from (\ref{GN11A1}) to (\ref{exponent11A2}) hold for $q$. From the definition of the norm in $X(t)$, we may conclude immediately the inequality (\ref{pt4.5}).
\par In the second step, let us prove the inequality (\ref{pt4.4}). Taking into consideration two elements $(u,v)$ and $(\bar{u},\bar{v})$ from $X(t)$ we have
$$N(u,v)(t,x)- N(\bar{u},\bar{v})(t,x)= \big(u^{nl}(t,x)- \bar{u}^{nl}(t,x), v^{nl}(t,x)- \bar{v}^{nl}(t,x)\big). $$
We use, on the one hand, the $(L^1 \cap L^2)- L^2$ estimates from Proposition \ref{md2.1} for $u^{nl}- \bar{u}^{nl}$ and $v^{nl}- \bar{v}^{nl}$. On the other hand, for $\partial^j_t|D|^{k\sigma_1}\big(u^{nl}- \bar{u}^{nl}\big)$ and $\partial^j_t|D|^{k\sigma_2}\big(v^{nl}- \bar{v}^{nl}\big)$, with $(j,k)=(0,1)$ or $(1,0)$, we apply $(L^1 \cap L^2)- L^2$ estimates if $\tau \in [0,t/2]$ and the $L^2-L^2$ estimates if $\tau \in [t/2,t]$ from Proposition \ref{md2.1}. Hence, we derive the following estimates for $(j,k)=(0,1)$ or $(1,0)$:
\begin{align*}
\big\|\big(u^{nl}- \bar{u}^{nl}\big)(t,\cdot)\big\|_{L^2} &\lesssim \int_0^t (1+t-\tau)^{-\frac{n}{4\sigma_1}}\big\||v(\tau,\cdot)|^p- \bar{v}(\tau,\cdot)|^p\big\|_{L^1 \cap L^2}d\tau, \\
\big\|\partial^j_t|D|^{k\sigma_1}\big(u^{nl}- \bar{u}^{nl}\big)(t,\cdot)\big\|_{L^2} &\lesssim \int_0^{t/2}(1+t-\tau)^{-\frac{n}{4\sigma_1}- \frac{k}{2}-j}\big\||v(\tau,\cdot)|^p- \bar{v}(\tau,\cdot)|^p\big\|_{L^1 \cap L^2}d\tau \\
&\qquad + \int_{t/2}^t (1+t-\tau)^{- \frac{k}{2}-j}\big\||v(\tau,\cdot)|^p- \bar{v}(\tau,\cdot)|^p\big\|_{L^2}d\tau,
\end{align*}
and
\begin{align*}
\big\|\big(v^{nl}- \bar{v}^{nl}\big)(t,\cdot)\big\|_{L^2} &\lesssim \int_0^t (1+t-\tau)^{-\frac{n}{4\sigma_2}}\big\||u(\tau,\cdot)|^q- \bar{u}(\tau,\cdot)|^q\big\|_{L^1 \cap L^2}d\tau, \\
\big\|\partial^j_t|D|^{k\sigma_2}\big(v^{nl}- \bar{v}^{nl}\big)(t,\cdot)\big\|_{L^2} &\lesssim \int_0^{t/2}(1+t-\tau)^{-\frac{n}{4\sigma_2}- \frac{k}{2}-j}\big\||u(\tau,\cdot)|^q- \bar{u}(\tau,\cdot)|^q\big\|_{L^1 \cap L^2}d\tau, \\
&\qquad + \int_{t/2}^t (1+t-\tau)^{- \frac{k}{2}-j}\big\||u(\tau,\cdot)|^q- \bar{u}(\tau,\cdot)|^q\big\|_{L^2}d\tau.
\end{align*}
After employing H\"{o}lder's inequality, we arrive at
\begin{align*}
\big\||v(\tau,\cdot)|^p- |\bar{v}(\tau,\cdot)|^p\big\|_{L^1} &\lesssim \|v(\tau,\cdot)- \bar{v}(\tau,\cdot)\|_{L^p} \big(\|v(\tau,\cdot)\|^{p-1}_{L^p}+ \|\bar{v}(\tau,\cdot)\|^{p-1}_{L^p}\big), \\
\big\||v(\tau,\cdot)|^p- |\bar{v}(\tau,\cdot)|^p\big\|_{L^2} &\lesssim \|v(\tau,\cdot)- \bar{v}(\tau,\cdot)\|_{L^{2p}} \big(\|v(\tau,\cdot)\|^{p-1}_{L^{2p}}+ \|\bar{v}(\tau,\cdot)\|^{p-1}_{L^{2p}}\big), \\
\big\||u(\tau,\cdot)|^q- |\bar{u}(\tau,\cdot)|^q\big\|_{L^1} &\lesssim \|u(\tau,\cdot)- \bar{u}(\tau,\cdot)\|_{L^q} \big(\|u(\tau,\cdot)\|^{q-1}_{L^q}+ \|\bar{u}(\tau,\cdot)\|^{q-1}_{L^q}\big), \\
\big\||u(\tau,\cdot)|^q- |\bar{u}(\tau,\cdot)|^q\big\|_{L^2} &\lesssim \|u(\tau,\cdot)- \bar{u}(\tau,\cdot)\|_{L^{2q}} \big(\|u(\tau,\cdot)\|^{q-1}_{L^{2q}}+ \|\bar{u}(\tau,\cdot)\|^{q-1}_{L^{2q}}\big).
\end{align*}
By using analogous arguments as we proceeded to prove (\ref{pt4.5}), we apply the fractional Gagliardo-Nirenberg inequality from Proposition \ref{fractionalGagliardoNirenberg} to the terms
\begin{align*}
&\|v(\tau,\cdot)-\bar{v}(\tau,\cdot)\|_{L^{\eta_1}},\quad \|u(\tau,\cdot)-\bar{u}(\tau,\cdot)\|_{L^{\eta_2}}, \\
&\|v(\tau,\cdot)\|_{L^{\eta_1}},\quad \|\bar{v}(\tau,\cdot)\|_{L^{\eta_1}},\quad \|u(\tau,\cdot)\|_{L^{\eta_2}},\quad \|\bar{u}(\tau,\cdot)\|_{L^{\eta_2}},
\end{align*}
with $\eta_1=p$ or $\eta_1=2p$, and $\eta_2=q$ or $\eta_2=2q$ to conclude the inequality (\ref{pt4.4}). Summarizing, the proof of Theorem \ref{dl1.1} is completed.

\subsection{Proof of Theorem \ref{dl1.2}}
We follow the proof of Theorem \ref{dl1.1} with minor modifications in the steps of our proof. We also introduce both spaces for the data and the solutions as in Theorem \ref{dl1.1}, where the weights (\ref{pt4.1}) and (\ref{pt4.2}) are modified in the following way:
\begin{align*}
&f_1(\tau)= (1+\tau)^{-\frac{n}{4\sigma_1}},\quad f_2(\tau)= (1+\tau)^{-\frac{n}{4\sigma_1}-\frac{1}{2}}, \quad f_3(\tau)=(1+\tau)^{-\frac{n}{4\sigma_1}-1}, \\
&g_1(\tau)= (1+\tau)^{-\frac{n}{4\sigma_2}+ \e(q,\sigma_1)},\quad g_2(\tau)= (1+\tau)^{-\frac{n}{4\sigma_2}-\frac{1}{2}+ \e(q,\sigma_1)}, \quad g_3(\tau)=(1+\tau)^{-\frac{n}{4\sigma_2}-1+ \e(q,\sigma_1)}.
\end{align*}
Then, repeating some steps of the proofs we did in Theorem \ref{dl1.1} we may complete the proof of Theorem \ref{dl1.2}.

\section{Proof of blow-up result} \label{Sec.Blow-up}
Our goal of this section is to find the critical exponents in our main results. The proof of blow-up result is based on a contradiction argument by using the test function method (see, for example, \cite{NishiharaWakasugi,SunWang}). When $\sigma$ is an integer number, we are going to apply standard test functions, i.e. test functions with compact support. However, it seems difficult to directly apply this strategy to the fractional Laplacian operators $(-\Delta)^\sigma$ as well-known non-local operators when $\sigma$ is a fractional number. For this reason, the application of a modified test function method from Section \ref{Sec.Pre} comes into play. We shall divide the proof of Theorem \ref{dl1.3} into two cases as follows.

\subsection{The case that $\sigma$ is integer} \label{Sec.4.1}
\begin{proof}
At first, we introduce the test functions $\eta= \eta(t)$ and $\varphi= \varphi(x)$ having the following properties:
\begin{align}
&1.\quad \eta \in \mathcal{C}_0^\ity([0,\ity)) \text{ and }
\eta(t)=\begin{cases}
1 &\quad \text{ for }0 \le t \le \frac{1}{2}, \\
\text{decreasing } &\quad \text{ for }\frac{1}{2} \le t \le 1, \\
0 &\quad \text{ for }t \ge 1,
\end{cases} \nonumber \\
&2.\quad \varphi \in \mathcal{C}_0^\ity(\R^n) \text{ and }
\varphi(x)= \begin{cases}
1 &\text{ for } |x|\le 1/2, \\
0 &\text{ for }|x|\ge 1,
\end{cases} & \nonumber \\
&3.\quad \eta^{-\frac{\kappa'}{\kappa}}(t)\big(|\eta'(t)|^{\kappa'}+|\eta''(t)|^{\kappa'}\big) \le C \quad \text{ for any } t \in \Big[\frac{1}{2},1\Big], \label{t13.1.1} \\
&\qquad \text{ and }\quad \varphi^{-\frac{\kappa'}{\kappa}}(x) |\Delta^{\sigma}\varphi(x)|^{\kappa'} \le C \quad \text{ for any } |x| \in \Big[\frac{1}{2},1\Big], \label{t13.1.2}
\end{align}
with $\kappa= p$ or $\kappa= q$, where $\kappa'$ is the conjugate of $\kappa$ and $C$ is a suitable positive constant. In addition, we suppose that $\varphi= \varphi(|x|)$ is a radial function satisfying $\varphi(|x|) \le \varphi(|y|)$ for any $|x|\ge |y|$.
\par Let $R$ be a large parameter in $[0,\ity)$. We define the following test function:
$$ \phi_R(t,x):= \eta_R(t) \varphi_R(x), $$
where $\eta_R(t):= \eta(R^{-2\sigma}t)$ and $\varphi_R(x):= \varphi(R^{-1}x)$. Now we define the functionals
\begin{align*}
&I_R:= \int_0^{\ity}\int_{\R^n}|v(t,x)|^p \phi_R(t,x) dxdt= \int_{Q_R}|v(t,x)|^p \phi_R(t,x) d(x,t),  \\ 
&J_R:= \int_0^{\ity}\int_{\R^n}|u(t,x)|^q \phi_R(t,x) dxdt= \int_{Q_R}|u(t,x)|^q \phi_R(t,x) d(x,t),
\end{align*}
where $$Q_R:= \big[0,R^{2\sigma}\big] \times B_R \quad \text{ with }\quad B_R:= \big\{x\in \R^n: |x|\le R \big\}. $$
Let us assume that $(u,v)= \big(u(t,x),v(t,x)\big)$ is a global (in time) Sobolev solution to (\ref{pt1.1}). After multiplying the first equation to (\ref{pt1.1}) by $\phi_R=\phi_R(t,x)$, we carry out partial integration to obtain
\begin{align} 
&I_R+ \int_{B_R} u_1(x) \varphi_R(x)dx \nonumber \\
&\qquad= \int_{Q_R}u(t,x) \Big(\eta''_R(t) \varphi_R(x)- \eta'_R(t)  \varphi_R(x)+ \eta_R(t) (-\Delta)^{\sigma}\varphi_R(x)\Big)d(x,t). \label{t13.1.3}
\end{align}
Employing H\"{o}lder's inequality with $\frac{1}{q}+\frac{1}{q'}=1$ we can proceed as follows:
\begin{align*}
&\int_{Q_R} |u(t,x)|\, \big|\eta''_R(t) \varphi_R(x)\big| d(x,t) \\
&\qquad \le \Big(\int_{Q_R} \Big|u(t,x)\phi^{\frac{1}{q}}_R(t,x)\Big|^q d(x,t)\Big)^{\frac{1}{q}} \Big(\int_{Q_R} \Big|\phi^{-\frac{1}{q}}_R(t,x) \eta''_R(t) \varphi_R(x)\Big|^{q'} d(x,t)\Big)^{\frac{1}{q'}} \\
&\qquad \le J_R^{\frac{1}{q}} \Big( \int_{Q_R}\eta_R^{-\frac{q'}{q}}(t) \big|\eta''_R(t)\big|^{q'} \varphi_R(x) d(x,t)\Big)^{\frac{1}{q'}}.
\end{align*}
After performing change of variables $\tilde{t}:= R^{-2\sigma}t$ and $\tilde{x}:= R^{-1}x$, we derive
\begin{equation} \label{t13.1.4}
\int_{Q_R} |u(t,x)|\, \big|\eta''_R(t) \varphi_R(x)\big| d(x,t) \lesssim J_R^{\frac{1}{q}}\, R^{-4\sigma+ \frac{n+2\sigma}{q'}},
\end{equation}
where we used the relation $\eta''_R(t)= R^{-4\sigma}\eta''(\tilde{t})$ and the assumption (\ref{t13.1.1}). Analogously, we may arrive at the following estimates:
\begin{align}
&\int_{Q_R}|u(t,x)|\, \big|\eta'_R(t) \varphi_R(x)\big| d(x,t) \lesssim J_R^{\frac{1}{q}}\, R^{-2\sigma+ \frac{n+2\sigma}{q'}}, \label{t13.1.5} \\ 
&\int_{Q_R}|u(t,x)|\, \big|\eta_R(t) (-\Delta)^{\sigma}\varphi_R(x)\big| d(x,t) \lesssim J_R^{\frac{1}{q}}\, R^{-2\sigma+ \frac{n+2\sigma}{q'}}, \label{t13.1.6}
\end{align}
where we used the relations $\eta'_R(t)= R^{-2\sigma}\eta'(\tilde{t})$, $(-\Delta)^{\sigma}\varphi_R(x)= R^{-2\sigma}(-\Delta)^{\sigma}\varphi(\tilde{x})$ (since $\sigma$ is an integer number), and the assumptions (\ref{t13.1.1}) and (\ref{t13.1.2}). Due to the assumption (\ref{optimal13.1}), there exists a sufficiently large constant $R_0>0$ so that it holds
\begin{equation} \label{t13.1.7}
\int_{B_R} u_1(x) \varphi_R(x)dx >0
\end{equation}
for any $R>R_0$. As a result, combining the estimates from (\ref{t13.1.3}) to (\ref{t13.1.7}) we have proved that
$$ I_R \lesssim J_R^{\frac{1}{q}}\, R^{-2\sigma+ \frac{n+2\sigma}{q'}}. $$
In the same way we may conclude
$$ J_R \lesssim I_R^{\frac{1}{p}}\, R^{-2\sigma+ \frac{n+2\sigma}{p'}}. $$
Therefore, from the both above estimates we get
\begin{align}
I_R^{\frac{pq-1}{pq}} &\lesssim R^{(-2\sigma+ \frac{n+2\sigma}{p'})\frac{1}{q}- 2\sigma+ \frac{n+2\sigma}{q'}}=: R^{\gamma_1}, \label{t13.1.8} \\ 
J_R^{\frac{pq-1}{pq}} &\lesssim R^{(-2\sigma+ \frac{n+2\sigma}{q'})\frac{1}{p}- 2\sigma+ \frac{n+2\sigma}{p'}}=: R^{\gamma_2} \label{t13.1.9}.
\end{align}
Without loss of generality we can assume $q>p$. It follows immediately that the assumption (\ref{optimal13.2}) is equivalent to
$$ \frac{n}{2\sigma} \le \frac{1+q}{pq-1}, $$
that is, $\gamma_2 \le 0$. For this reason, we shall split our consideration into two subcases.
\par In the subcritical case $\frac{n}{2\sigma}< \frac{1+q}{pq-1}$, i.e. $\gamma_2 <0$, passing $R \to \ity$ in (\ref{t13.1.9}) we obtain
$$ \int_0^{\ity}\int_{\R^n}|u(t,x)|^q dxdt= 0. $$
This implies $u \equiv 0$, which is a contradiction to the assumption (\ref{optimal13.1}). Therefore, there is no global (in time) Sobolev solution to (\ref{pt1.1}) in the subcritical case.
\par Let us consider the critical case $\frac{n}{2\sigma}= \frac{1+q}{pq-1}$, i.e. $\gamma_2=0$. From (\ref{t13.1.9}) there exists a positive constant $C_0$ such that
$$ J_R= \int_{Q_R}|u(t,x)|^q \phi_R(t,x) d(x,t) \le C_0, $$
for a sufficiently large $R$. Thus, it follows that
\begin{equation} \label{t13.1.10}
\int_{\bar{Q}_R}|u(t,x)|^q \phi_R(t,x) d(x,t) \to 0 \quad \text{ as } R \to \ity,
\end{equation}
where we introduce the notations
$$ \bar{Q}_R:= Q_R \setminus \big(\big[0,R^{2\sigma}/2\big] \times B_{R/2}\big) \quad \text{ with }\quad B_{R/2}:= \big\{x\in \R^n: 0\le |x|\le R/2 \big\}. $$
Since $\partial_t^2 \phi_R(t,x)= \partial_t \phi_R(t,x)= (-\Delta)^{\sigma}\phi_R(t,x)=0$ in $\big(\R^1_+ \times \R^n\big) \setminus \bar{Q}_R$, we may repeat the steps of the proofs from (\ref{t13.1.3}) to (\ref{t13.1.6}) to conclude the following estimates:
\begin{align*}
I_R+ \int_{B_R} u_1(x) \varphi_R(x)dx&\lesssim \Big(\int_{\bar{Q}_R}|u(t,x)|^q \phi_R(t,x) d(x,t)\Big)^{\frac{1}{q}}\, R^{-2\sigma+ \frac{n+2\sigma}{q'}}, \\ 
J_R+ \int_{B_R} v_1(x) \varphi_R(x)dx&\lesssim \Big(\int_{\bar{Q}_R}|v(t,x)|^p \phi_R(t,x) d(x,t)\Big)^{\frac{1}{p}}\, R^{-2\sigma+ \frac{n+2\sigma}{p'}}.
\end{align*}
Because $\gamma_2= 0$, from the both above estimates and (\ref{t13.1.7}) we arrive at
\begin{equation} \label{t13.1.11}
J_R+ \int_{B_R} v_1(x) \varphi_R(x)dx \lesssim \Big(\int_{\bar{Q}_R}|u(t,x)|^q \phi_R(t,x) d(x,t)\Big)^{\frac{1}{pq}}.
\end{equation}
By using (\ref{t13.1.10}) we let $R \to \ity$ in (\ref{t13.1.11}) to derive
$$ \int_0^{\ity}\int_{\R^n}|u(t,x)|^q dxdt+ \int_{\R^n} v_1(x) dx= 0. $$
This is again a contradiction to the assumption (\ref{optimal13.1}), that is, there is no global (in time) Sobolev solution to (\ref{pt1.1}) in the critical case. Hence, this completes our proof.
\end{proof}

\subsection{The case that $\sigma$ is fractional} \label{Sec.4.2}
\begin{proof}
At first, we denote the constant $\theta:= \sigma-[\sigma]$. Since $\sigma$ is a fractional number, it follows that $\theta \in (0,1)$. We introduce, on the one hand, the function $\eta= \eta(t)$ satisfying the same properties as in Section \ref{Sec.4.1}. On the other hand, we introduce the function $\varphi=\varphi(|x|):=\big< x\big>^{-n-2\theta}$.
\par Let $R$ be a large parameter in $[0,\ity)$. We define the following test function:
$$ \phi_R(t,x):= \eta_R(t) \varphi_R(x), $$
where $\eta_R(t):= \eta(R^{-2\sigma}t)$ and $\varphi_R(x):= \varphi(R^{-1}x)$. We define the funtionals
\begin{align*}
&I_R:= \int_0^{\ity}\int_{\R^n}|v(t,x)|^p \phi_R(t,x)\,dxdt= \int_0^{R^{2\sigma}}\int_{\R^n}|v(t,x)|^p \phi_R(t,x)\,dxdt, \\
&J_R:=  \int_0^{\ity}\int_{\R^n}|u(t,x)|^q \phi_R(t,x)\,dxdt= \int_0^{R^{2\sigma}}\int_{\R^n}|u(t,x)|^q \phi_R(t,x)\,dxdt,
\end{align*}
and
$$ I_{R,t}:= \int_{\frac{R^\alpha}{2}}^{R^{2\sigma}}\int_{\R^n}|v(t,x)|^p \phi_R(t,x)\,dxdt, \qquad J_{R,t}:= \int_{\frac{R^{2\sigma}}{2}}^{R^{2\sigma}}\int_{\R^n}|u(t,x)|^q \phi_R(t,x)\,dxdt. $$
Let us assume that $(u,v)= \big(u(t,x),v(t,x)\big)$ is a global (in time) Sobolev solution from $C\big([0,\infty),L^2\big) \times C\big([0,\infty),L^2\big)$ to (\ref{pt1.1}). After multiplying the first equation to (\ref{pt1.1}) by $\phi_R=\phi_R(t,x)$, we perform partial integration to obtain
\begin{align}
0\le I_R &= -\int_{\R^n} u_1(x)\varphi_R(x)\,dx + \int_{\frac{R^{2\sigma}}{2}}^{R^{2\sigma}}\int_{\R^n}u(t,x) \eta''_R(t) \varphi_R(x)\,dxdt \nonumber \\
&\quad + \int_0^{\ity}\int_{\R^n} \eta_R(t) \varphi_R(x)\, (-\Delta)^{\sigma} u(t,x)\,dxdt- \int_{\frac{R^{2\sigma}}{2}}^{R^{2\sigma}}\int_{\R^n} \eta'_R(t) \varphi_R(x)\, u(t,x)\,dxdt \nonumber \\
&=: -\int_{\R^n} u_1(x)\varphi_R(x)\,dx+ I_{1R}+ I_{2R}- I_{3R}. \label{t13.2.1}
\end{align}
Applying H\"{o}lder's inequality with $\frac{1}{q}+\frac{1}{q'}=1$ we may deal with as follows:
\begin{align*}
|I_{1R}| &\le \int_{\frac{R^{2\sigma}}{2}}^{R^{2\sigma}}\int_{\R^n} |u(t,x)|\, \big|\eta''_R(t)\big| \varphi_R(x) \, dxdt \\
&\lesssim \Big(\int_{\frac{R^{2\sigma}}{2}}^{R^{2\sigma}}\int_{\R^n} \Big|u(t,x)\phi^{\frac{1}{q}}_R(t,x)\Big|^p \,dxdt\Big)^{\frac{1}{q}} \Big(\int_{\frac{R^{2\sigma}}{2}}^{R^{2\sigma}}\int_{\R^n} \Big|\phi^{-\frac{1}{q}}_R(t,x) \eta''_R(t) \varphi_R(x)\Big|^{p'}\, dxdt\Big)^{\frac{1}{q'}} \\
&\lesssim J_{R,t}^{\frac{1}{q}}\, \Big(\int_{\frac{R^{2\sigma}}{2}}^{R^{\alpha}}\int_{\R^n} \eta_R^{-\frac{q'}{q}}(t) \big|\eta''_R(t)\big|^{q'} \varphi_R(x)\, dxdt\Big)^{\frac{1}{q'}}.
\end{align*}
By using the change of variables $\tilde{t}:= R^{-2\sigma}t$ and $\tilde{x}:= R^{-1}x$, we compute directly to give
\begin{equation}
|I_{1R}| \lesssim J_{R,t}^{\frac{1}{q}}\, R^{-4\sigma+ \frac{n+2\sigma}{q'}}\Big(\int_{\R^n} \big< \tilde{x}\big>^{-n-2\theta}\, d\tilde{x}\Big)^{\frac{1}{q'}}. \label{t13.2.2}
\end{equation}
Here we used $\eta''_R(t)= R^{-4\sigma}\eta''(\tilde{t})$ and the assumption (\ref{t13.1.1}). In an analogous way, we may conclude the following estimate for $I_{3R}$:
\begin{equation}
|I_{3R}| \lesssim J_{R,t}^{\frac{1}{q}}\, R^{-2\sigma+ \frac{n+2\sigma}{q'}}\Big(\int_{\R^n} \big< \tilde{x}\big>^{-n-2\theta}\, d\tilde{x}\Big)^{\frac{1}{q'}}. \label{t13.2.3}
\end{equation}
Now let us focus our attention to estimate $I_{2R}$. In the first step, since $\varphi_R \in H^{2\sigma}$ and $u \in C\big([0,\infty),L^2\big)$, we apply Lemma \ref{lemma2.3} to derive the following relations:
$$ \int_{\R^n} \varphi_R(x)\, (-\Delta)^{\sigma} u(t,x)\,dx= \int_{\R^n}|\xi|^{2\sigma} \widehat{\varphi}_R(\xi)\,\widehat{u}(t,\xi)\,d\xi= \int_{\R^n} u(t,x)\, (-\Delta)^{\sigma} \varphi_R(x)\,dx. $$
Therefore, we get
$$I_{2R}= \int_0^{\ity}\int_{\R^n} \eta_R(t) \varphi_R(x)\, (-\Delta)^{\sigma} u(t,x)\,dxdt= \int_0^{\ity}\int_{\R^n} \eta_R(t) u(t,x)\, (-\Delta)^{\sigma} \varphi_R(x) \,dxdt. $$
Employing H\"{o}lder's inequality again as we did to estimate $J_1$ gives
$$|I_{2R}|\le I_R^{\frac{1}{q}}\, \Big(\int_0^{R^{2\sigma}}\int_{\R^n} \eta_R(t) \varphi^{-\frac{q'}{q}}_R(x)\, \big|(-\Delta)^{\sigma} \varphi_R(x)\big|^{q'} \, dxdt\Big)^{\frac{1}{q'}}. $$
In the second step, to control the above integral we shall apply results from Lemmas \ref{lemma2.1} and \ref{lemma2.2} as the key tools. In particular, carrying out the change of variables $\tilde{x}:= R^{-1}x$ we get the following relation from Lemma \ref{lemma2.2}:
$$ (-\Delta)^{\sigma}\varphi_R(x)= R^{-2\sigma} (-\Delta)^{\sigma}(\varphi)(\tilde{x}). $$
For this reason, using the change of variables $\tilde{t}:= R^{-2\sigma}t$ leads to
\begin{align*}
|I_{2R}| &\lesssim J_R^{\frac{1}{q}}\, R^{-2\sigma+ \frac{n+2\sigma}{q'}}\Big(\int_0^{1}\int_{\R^n} \eta(\tilde{t}) \varphi^{-\frac{q'}{q}}(\tilde{x})\, \big|(-\Delta)^{\sigma} (\varphi)(\tilde{x})\big|^{q'}\, d\tilde{x}d\tilde{t}\Big)^{\frac{1}{q'}} \\
&\lesssim J_R^{\frac{1}{q}}\, R^{-2\sigma+ \frac{n+2\sigma}{q'}}\Big(\int_{\R^n} \varphi^{-\frac{q'}{q}}(\tilde{x})\, \big|(-\Delta)^{\sigma} (\varphi)(\tilde{x})\big|^{q'}\, d\tilde{x}\Big)^{\frac{1}{q'}}.
\end{align*}
After employing Lemma \ref{lemma2.1}, we deduce the following estimate:
\begin{equation}
|I_{2R}|\lesssim I_R^{\frac{1}{q}}\, R^{-2\sigma+ \frac{n+2\sigma}{q'}}\Big(\int_{\R^n} \big< \tilde{x}\big>^{-n-2\theta}\, d\tilde{x}\Big)^{\frac{1}{q'}}. \label{t13.2.4}
\end{equation}
Due to the assumption (\ref{optimal13.1}), there exists a sufficiently large constant $R_1> 0$ such that it holds
\begin{equation}
\int_{\R^n} u_1(x) \varphi_R(x)\, dx >0 \label{t13.2.5}
\end{equation}
for all $R > R_1$. Combining the estimates from (\ref{t13.2.1}) to (\ref{t13.2.5}) we may arrive at
\begin{align}
0< \int_{\R^n} u_1(x) \varphi_R(x)\, dx &\lesssim J_{R,t}^{\frac{1}{q}} \Big(R^{-4\sigma+ \frac{n+2\sigma}{q'}}+ R^{-2\sigma+ \frac{n+2\sigma}{q'}}\Big)+ J_R^{\frac{1}{q}}\, R^{-2\sigma+ \frac{n+2\sigma}{q'}}- I_R \nonumber \\
&\lesssim J_R^{\frac{1}{q}} R^{-2\sigma+ \frac{n+2\sigma}{q'}}- I_R \label{t13.2.6}
\end{align}
for all $R > R_1$. Similarly, we may conclude the following estimate for all $R > R_1$:
\begin{align}
0< \int_{\R^n} v_1(x) \varphi_R(x)\, dx &\lesssim I_{R,t}^{\frac{1}{p}} \Big(R^{-4\sigma+ \frac{n+2\sigma}{p'}}+ R^{-2\sigma+ \frac{n+2\sigma}{p'}}\Big)+ I_R^{\frac{1}{p}}\, R^{-2\sigma+ \frac{n+2\sigma}{p'}}- J_R \nonumber \\
&\lesssim I_R^{\frac{1}{p}} R^{-2\sigma+ \frac{n+2\sigma}{p'}}- J_R. \label{t13.2.7}
\end{align}
As a result, combining the estimates (\ref{t13.2.6}) and (\ref{t13.2.7}) we derive
\begin{align*}
I_R&\lesssim J_R^{\frac{1}{q}}R^{-2\sigma+ \frac{n+2\sigma}{q'}}, \\ 
J_R&\lesssim I_R^{\frac{1}{p}}R^{-2\sigma+ \frac{n+2\sigma}{p'}}.
\end{align*}
Thus, it follows immediately that
\begin{align}
I_R^{\frac{pq-1}{pq}} &\lesssim R^{- 2\sigma+ \frac{n+2\sigma}{q'}+ (-2\sigma+ \frac{n+2\sigma}{p'})\frac{1}{q}}=: R^{\gamma_1}, \label{t13.2.8} \\ 
J_R^{\frac{pq-1}{pq}} &\lesssim R^{- 2\sigma+ \frac{n+2\sigma}{p'}+ (-2\sigma+ \frac{n+2\sigma}{q'})\frac{1}{p}}=: R^{\gamma_2} \label{t13.2.9}.
\end{align}
Without loss of generality we can assume $q>p$. We can see that the assumption (\ref{optimal13.2}) is equivalent to $\gamma_2 \le 0$. For this reason, we shall split our consideration into two subcases.
\par In the subcritical case $\frac{n}{2\sigma}< \frac{1+q}{pq-1}$, i.e. $\gamma_2 <0$, we pass $R \to \ity$ in (\ref{t13.2.9}) to get
$$ J_R=  \int_0^{\ity}\int_{\R^n}|u(t,x)|^q \phi_R(t,x)\,dxdt= 0. $$
which follows $u \equiv 0$. This is a contradiction to the assumption (\ref{optimal13.1}), that is, there is no global (in time) Sobolev solution to (\ref{pt1.1}) in the subcritical case.
\par Let us turn to the critical case $\frac{1+q}{pq-1}= \frac{n}{2\sigma}$, i.e. $\gamma_2= 0$. First, we introduce the following constants:
\begin{align*}
&C_{0u}:= \int_{\R^n} u_1(x) \varphi_R(x) \quad \text{ and }\quad C_{0v}:= \int_{\R^n} v_1(x) \varphi_R(x), \\ 
&C_{p'}:= \Big(\int_{\R^n} \big< \tilde{x}\big>^{-n-2\theta}\, d\tilde{x}\Big)^{\frac{1}{p'}} \quad \text{ and }\quad C_{q'}:= \Big(\int_{\R^n} \big< \tilde{x}\big>^{-n-2\theta}\, d\tilde{x}\Big)^{\frac{1}{q'}}.
\end{align*}
Then, repeating some arguments as we did in the subcritical case we may arrive at the following estimates:
\begin{align*}
&0< I_R+ C_{0u} \le C_{q'} J_R^{\frac{1}{q}}R^{-2\sigma+ \frac{n+2\sigma}{q'}}, \\ 
&0< J_R+ C_{0v} \le C_{p'} I_R^{\frac{1}{p}}R^{-2\sigma+ \frac{n+2\sigma}{p'}}.
\end{align*}
As a result, we deduce that
\begin{equation}
J_R+ C_{0v} \le C_{p'}C^{\frac{1}{p}}_{q'} J_R^{\frac{1}{pq}}R^{\gamma_2}= C_{p'}C^{\frac{1}{p}}_{q'} J_R^{\frac{1}{pq}}. \label{t13.2.10}
\end{equation}
It is clear that the estimate (\ref{t13.2.10}) follows $J_R \le C_{p'}C^{\frac{1}{p}}_{q'} J_R^{\frac{1}{pq}}$ and $C_{0v} \le C_{p'}C^{\frac{1}{p}}_{q'} J_R^{\frac{1}{pq}}$. For this reason, we derive
\begin{equation}
J_R \le C_0, \label{t13.2.11}
\end{equation}
where $C_0:= \Big(C_{p'}C^{\frac{1}{p}}_{q'}\Big)^{\frac{pq}{pq-1}}$ is a positive constant, and
\begin{equation}
J_R \ge \left(\frac{C_{0v}}{C_{p'}C^{\frac{1}{p}}_{q'}}\right)^{pq}. \label{t13.2.12}
\end{equation}
After plugging (\ref{t13.2.12}) into the left-hand side of (\ref{t13.2.10}), we compute straightforwardly to obtain
$$ J_R \ge \frac{(C_{0v})^{(pq)^2}}{\Big(C_{p'}C^{\frac{1}{p}}_{q'}\Big)^{pq+(pq)^2}}. $$
For any integer $j\ge 1$, using an iteration argument gives
\begin{equation}
J_R \ge \frac{C_{0v}^{(pq)^j}}{\Big(C_{p'}C^{\frac{1}{p}}_{q'}\Big)^{pq+(pq)^2+\cdots+(pq)^j}}= \frac{C_{0v}^{(pq)^j}}{\Big(C_{p'}C^{\frac{1}{p}}_{q'}\Big)^{\frac{(pq)^{j+1}-pq}{pq-1}}}= \Big(C_{p'}C^{\frac{1}{p}}_{q'}\Big)^{\frac{pq}{pq-1}}\left(\frac{C_{0v}}{\Big(C_{p'}C^{\frac{1}{p}}_{q'}\Big)^{\frac{pq}{pq-1}}}\right)^{(pq)^j}. \label{t13.2.13}
\end{equation}
Let us now choose the constant
$$ \epsilon_2= \int_{\R^n} \big< \tilde{x}\big>^{-n-2\theta}\, d\tilde{x} $$
in the assumption (\ref{optimal13.1}). Then, there exists a sufficiently large constant $R_2> 0$ so that
$$ \int_{\R^n} v_1(x) \varphi_R(x)\, dx > \epsilon_2 $$
for all $R > R_2$. It is obvious that this is equivalent to
$$ C_{0v}> \int_{\R^n} \big< \tilde{x}\big>^{-n-2\theta}\, d\tilde{x}= \Big(C_{p'}C^{\frac{1}{p}}_{q'}\Big)^{\frac{pq}{pq-1}}, \quad \text{that is, }\quad \frac{C_{0v}}{\Big(C_{p'}C^{\frac{1}{p}}_{q'}\Big)^{\frac{pq}{pq-1}}}> 1. $$
Letting $j \to \ity$ in (\ref{t13.2.13}) we may conclude $J_R \to \ity$, which is a contradiction to (\ref{t13.2.11}). Therefore, there is no global (in time) Sobolev solution to (\ref{pt1.1}) in the critical case. Summarizing, the proof of Theorem \ref{dl1.3} is completed.
\end{proof}

\section*{Acknowledgment}
The PhD study of MSc. T.A. Dao is supported by Vietnamese Government's Scholarship (Grant number: 2015/911). The author would like to thank sincerely to Prof. Michael Reissig for valuable discussions and Institute of Applied Analysis for their hospitality. The author is grateful to the referee for his careful reading of the manuscript and for helpful comments. \medskip

\section*{Appendix}

\begin{md}[Fractional Gagliardo-Nirenberg inequality] \label{fractionalGagliardoNirenberg}
Let $1<p,\, p_0,\, p_1<\infty$, $\sigma >0$ and $s\in [0,\sigma)$. Then, it holds the following fractional Gagliardo-Nirenberg inequality for all $u\in L^{p_0} \cap \dot{H}^\sigma_{p_1}$:
$$ \|u\|_{\dot{H}^{s}_p}\lesssim \|u\|_{L^{p_0}}^{1-\theta}\,\, \|u\|_{\dot{H}^{\sigma}_{p_1}}^\theta, $$
where $\theta=\theta_{s,\sigma}(p,p_0,p_1)=\frac{\frac{1}{p_0}-\frac{1}{p}+\frac{s}{n}}{\frac{1}{p_0}-\frac{1}{p_1}+\frac{\sigma}{n}}$ and $\frac{s}{\sigma}\leq \theta\leq 1$.
\end{md}
For the proof one can see \cite{Ozawa}.

\bbd \label{UsefulLemma}
Let $s\in (0,1)$. Let $\psi$ be a smooth function satisfying $\partial_x^2 \psi\in L^\ity$. For any $R>0$, let $\psi_R$ be a function defined by
$$ \psi_R(x):= \psi\big(R^{-1} x\big)\quad \text{ for all }x \in \R^n. $$
Then, $(-\Delta)^s (\psi_R)$ satisfies the following scaling properties for all $x \in \R^n$:
$$(-\Delta)^s (\psi_R)(x)= R^{-2s}\big((-\Delta)^s \psi \big)\big(R^{-1} x\big). $$
\ebd
This result can be found in \cite{DaoReissig}.


\end{document}